\documentclass[letterpaper, 10pt, conference]{ieeeconf}
\IEEEoverridecommandlockouts 
\usepackage{amsmath,amssymb,url,subfigure,graphicx,theorem}
\usepackage{graphicx,subfigure,hyperref}
\usepackage{color,comment}
\usepackage{siunitx}
\usepackage{tikz}

\newcommand{\parenth}[1]{\ensuremath{\left( #1 \right)}}
\newcommand{\pair}[1]{\ensuremath{\langle #1 \rangle}}

\newcommand{\refeqn}[1]{(\ref{eqn:#1})}

\newcommand{\trs}[1]{\mathrm{tr}\ensuremath{[#1]}}

\newcommand{\SO}{\operatorname{SO}(3)}

\newcommand{\so}{\ensuremath{\mathfrak{so}(3)}}

\renewcommand{\Re}{\ensuremath{\mathbb{R}}}

\newcommand{\D}{\ensuremath{\mathbf{D}}}

\newcommand{\Sph}{\ensuremath{\mathsf{S}}}

\graphicspath{{./Figs/}}

\newtheorem{lemma}{Lemma}
\newtheorem{assumption}{Assumption}
\newtheorem{theorem}{Theorem}

\title{\Large\bf Attitude Observer on $\SO$ with Time-Varying Reference Directions}
\author{Kanishke Gamagedara, Taeyoung Lee and Dong Eui Chang%
\thanks{Kanishke Gamagedara and T. Lee are with the Department of Mechanical and Aerospace Engineering,  George Washington University, Washington DC, USA.
{\tt\small  \{kanishkegb,tylee\}@gwu.edu}}%
\thanks{D.E. Chang is with the School of Electrical Engineering, Korea Advanced Institute of Science and Technology (KAIST), Daejeon, South Korea.
        {\tt\small dechang@kaist.ac.kr}}%
    \thanks{This research been supported in part by NSF under the grants CMMI-1760928 and CNS-1837382, by the KUSTAR-KAIST Institute, KAIST, Korea and again by KAIST under the grants G04170001 and N11180231.}}

\begin{document}
\allowdisplaybreaks

\maketitle

\begin{abstract}
    This paper introduces an advanced Lyapunov stability analysis for an attitude observer that has been developed on the special orthogonal group. 
    In particular, when the attitude observer is constructed based on multiple direction measurements toward known reference points, a local exponential stability has been established by linearization, under the assumption that those reference points are fixed in the inertial frame. 
    Several modifications have been proposed to deal with reference directions changing over time. 
    Here, we present an alternative Lyapunov analysis to show that the attitude observer still exhibits exponential stability for time-varying reference directions, under the assumption that the observer gain is sufficiently large relative to the rate of change of the reference directions. 
    These are illustrated by a numerical example, followed by experimental results with visual marker detection in an indoor space. 
\end{abstract}

\section{Introduction}

There is a myriad of approaches in the development of attitude observers~\cite{CraMarJGCD07}, especially in terms of quaternions. 
However, it is well known that  there is an ambiguity in representing the attitude with quaternions, as the three-sphere, or the space of unit vectors in $\Re^4$ double covers the configuration space of the attitude defined as the special orthogonal group. 
When constructing an attitude control system or an attitude observer in terms of quaternions, an exogenous system is required to represent the attitude in a consistent manner, and otherwise,  undesired phenomena such as unwinding may appear~\cite{BhaBerSCL00}. 

To avoid these issues, attitude observers have been constructed directly on the special orthogonal group.
In particular, reference \cite{MahHamITAC08} presents a set of attitude observers comparable to nonlinear complementary filters under various assumptions. 
These attitude observers can be categorized by the following criteria: whether the attitude measurement is a set of direction measurements or complete attitudes; for the former, whether the reference direction is fixed in the inertial frame or not; whether a gyro bias is considered or not.

This paper focuses on the specific case when the measurement for attitude determination is given by multiple direction measurements toward known reference directions, which may vary over time or depending on the location of the vehicle. 
This is particularly useful for landmark based observers or indoor applications where the magnetic field is corrupted. 
While such cases can be addressed in a stochastic fashion with Bayesian framework~\cite{LeeITAC18}, we consider deterministic observers in this paper. 

For time-varying reference directions, stability analysis for attitude observers is presented in~\cite{TruMahITAC12}, with sufficient conditions for persistent-excitation that is also applicable to single direction measurements.
This work is based on assumption that there is no gyro bias. 
Later in~\cite{GriFosITAC12}, the attitude observer presented in~\cite{MahHamITAC08} is modified with a projection operation to deal with time-varying reference directions and a gyro bias concurrently. 

This paper presents an alternative stability analysis for the attitude observer of~\cite{MahHamITAC08}, and we show that it can handle time-varying reference directions without need for additional projection operation.
While this results in a reduced region of attraction compared with the projection based approach of~\cite{GriFosITAC12}, the presented stability analysis on the special orthogonal group can be utilized in the developments of hybrid attitude observers or non-memoryless observers~\cite{LeeChaAJDSMC18} to achieve global attractivity in time-varying reference directions. 

Next, the presented stability analysis is verified under indoor attitude estimation experiments, where the reference directions are comprised of the direction of gravity measured by an acceleration, and visual landmarks. 
More specifically, two feature points are placed in the environments, and they are observed by a low-cost camera.
Then, the video images are process via the OpenCV library to construct the line of sight measurements represented in the body-fixed frame. 
As the object undergoes both translations and rotations, the reference direction changes over time depending on the relative position toward feature points. 
The estimated attitudes are compared against the output of an external motion capture system with higher accuracy for validation. 

In short, the main contribution of this paper is an alternative, advanced stability analysis on the special orthogonal group for an attitude observer with time-varying reference directions and a gyro bias, and experimental implementation in attitude estimation with visual landmarks. 

\section{Problem Formulation}\label{sec:PF}

\subsection{Mathematical Preliminaries}

The inner product $\langle A, B\rangle $ of two matrices or vectors $A$ and $B$ of the same size denotes the usual Euclidean inner product, i.e., $\langle A, B\rangle = \operatorname{tr}(A^TB)$. 
The norm  $\| A\|$ for a matrix or vector $A$ denotes the Euclidean norm or the Frobenius norm, i.e., $\|A\|^2 = \langle A, A\rangle = \operatorname{tr}(A^TA)$. 
The minimum eigenvalue of a symmetric matrix $A$ is denoted by $\lambda_{\rm min} (A)$ and the maximum eigenvalue by $\lambda_{\rm max} (A)$.

Consider the attitude dynamics of a rigid body. 
We define the body-fixed frame and the inertial  reference frame. The attitude dynamics evolve on 
\[
\SO=\{ R\in\Re^{3\times 3}\,|\, R^TR=I_{3\times 3},\, \mathrm{det}[R]=1\},
\]
where the rotation matrix $R\in\SO$ corresponds to the linear transformation of the representation of a vector from the body-fixed frame to the inertial frame. 
For any $R, R_1, R_2 \in \SO$,
\begin{equation}\label{eqn:R:isometry:0}
 \langle RR_1, RR_2\rangle = \langle R_1, R_2 \rangle =\langle R_1R, R_2 R\rangle
\end{equation}
and
\begin{equation}\label{eqn:R_isometry}
\|RR_1 - RR_2\| = \| R_1  - R_2 \| = \|R_1R - R_2R\|.
\end{equation}
For any $R\in\SO$, there exist $\theta \in [0,2\pi]$ and a unit vector $v\in \Sph^2=\{q\in\Re^3\,|\, \|q\|=1\}$ so that
\[
R=\exp(\theta \hat v) =
I_{3\times 3} +\sin\theta\hat v + (1-\cos\theta)\hat v^2,
\] 
where the hat map $\wedge:\Re^3\rightarrow\so$ is defined such that $\hat x y = x\times y$ and $\hat x^T=-\hat x$ for any $x,y\in\Re^3$. 
The inverse of the hat map is denoted by the vee map $\vee:\so\rightarrow\Re^3$.

The following property is utilized in the subsequent development of the attitude observer. 
\begin{lemma}\label{lemma:IQ}
For any $Q \in \SO$ and $G=G^T \in \mathbb R^{3\times 3}$, 
\begin{equation}\label{eqn:IQ}
\langle G(I_{3\times 3}-Q), (I_{3\times 3}-Q)\rangle = 2(\trs{G} - \langle v_Q, Gv_Q\rangle)(1-\cos\theta) 
\end{equation}
where $\cos\theta = (\trs{Q} - 1)/2$ and $v_Q \in \Sph^2$ such that $Q = \exp(\theta \widehat{v_Q})$.
\end{lemma}
Throughout this paper, all of the proofs are relegated to Appendix.

\subsection{Attitude Observer Design Problem}

The attitude kinematics equation is given by
\begin{equation}
\dot R = R\hat\Omega,\label{eqn:R_dot}
\end{equation}
where $\Omega\in\Re^3$ is the angular velocity of the rigid body resolved in the body-fixed frame. 

There is an angular velocity sensor that measures the angular velocity of the rigid body up to a fixed bias. 
In other words, the measured angular velocity $\Omega_z\in\Re^3$ is given by
\begin{equation}
\Omega_z = \Omega + \gamma, \label{eqn:Wz}
\end{equation}
with a fixed bias $\gamma\in\Re^3$. 

It is assumed that the angular velocity is bounded.
\begin{assumption}
There is a positive constant $B_\Omega$ satisfying
\begin{equation}\label{eqn:BW}
\|\Omega(t)\| \leq B_\Omega,
\end{equation}
for any $t\geq 0$. 
\end{assumption}
For the determination of  attitude, suppose that there are $n$ distinctive objects, and the direction toward each of those objects is prescribed with respect to the inertial frame. 
Specifically, the direction to the $i$-th object in the inertial frame is given by the unit-vector $s_i(t)\in\Sph^2$, and it is assumed that $s_i(t)$ is available as a function of time $t$. 
For some positive weighting parameters $w_i$'s, define a time-varying, symmetric matrix $G(t)\in\Re^{3\times 3}$ as
\begin{equation}
    G(t)= \sum_{i=1}^n w_i s_i(t) s_i(t)^T.\label{eqn:G}
\end{equation}

By definition, the matrix $G(t)$ is positive-semidefinite always, and consequently, all of the eigenvalues of $G(t)$ are non-negative. 
Here we assume that the second largest eigenvalue of $G(t)$ is strictly positive as follows. 

\begin{assumption}
    Let $\lambda_2(G(t))\in\Re$ be the second largest eigenvalue of $G(t)$. 
    There is a positive constant $c$ such that
    \begin{equation}
        c \leq \lambda_2(G(t))\label{eqn:lambda_G}
    \end{equation}
for all $t\geq 0$.
\end{assumption}
This implies that $\mathrm{rank}(G(t)) \geq 2$ for all $t\geq0$. 
This is to ensure that there are at at least two non-parallel reference directions available, so that the attitude can be completely determined by direction measurements. 
We further assume that the rate of change of reference directions is bounded.
\begin{assumption}
    There is a positive constant $d>0$ such that
    \begin{align}
        \| \dot G(t)\| \leq d,\label{eqn:d}
    \end{align}
    for all $t\geq 0$.
\end{assumption}

There are sensors attached to the rigid body that can measure each direction of $s_i$. 
The sensor measurement to $s_i$ is given by $b_i\in\Sph^2$, and it is represented with respect to the body-fixed frame. 
Therefore,
\begin{equation}\label{eqn:bi}
b_i = R^T s_i,
\end{equation}
for $i\in\{1,\ldots n\}$. 

We wish to design an attitude observer to determine the attitude and the gyro bias. 
It is required that the observer is expressed in terms of the measurements of the directions and the angular velocity, without need for constructing the attitude directly from the direction measurements at every time instance.

\section{Attitude Observer on $\SO$}\label{sec:ASG}

\subsection{Error Variables}

Let $\bar R\in\SO$ and $\bar\gamma\in\Re^3$ be the estimated values of $R$ and $\gamma$, respectively. 
Define the estimation error variables $E_R\in\Re^{3\times 3}$ and $ e_\gamma \in\Re^3$ as
\begin{align}
E_R &= R- \bar R,\label{eqn:ER}\\
e_\gamma & = \gamma-\bar\gamma.\label{eqn:tilde_gamma}
\end{align}

Since the estimation error for the $i$-th reference direction is given by $
\|b_i-\bar R^T s_i\|^2= \|R^Ts_i-\bar R^T s_i\|^2= \|E_R^Ts_i\|^2 = \langle E_R^Ts_i, E_R^Ts_i \rangle =  \langle s_is_i^TE_R, E_R\rangle$, 
 it is natural to introduce the attitude estimation error function as 
\begin{align}\label{eqn:Psi}
\Psi(\bar R,t) 
& = \sum_{i=1}^n \frac{1}{2} w_i\langle s_is_i^TE_R, E_R\rangle
 = \frac{1}{2}\langle G E_R, E_R\rangle .
\end{align}

Several properties of the error variables are listed as follows. 

\begin{lemma}\label{lemma:important}
\begin{enumerate}\renewcommand{\labelenumi}{(\roman{enumi})}
\item The derivative of the error function $\Psi(\bar R,t)$ with respect to $\bar R$ along $\delta \bar R =\bar R\hat{\bar\eta}$ for $\bar\eta\in\Re^3$ is 
\begin{align}
\D_{\bar R} \Psi(\bar R,t) \cdot \delta \bar R  = \bar\eta \cdot e_R, \label{eqn:DPsi}
\end{align}
where the attitude error vector $e_R\in\Re^3$ is defined as
\begin{align}
e_R  & = (R^T G \bar R - \bar R^T G R)^\vee 
= \sum_{i=1}^n w_i \bar R^T s_i \times b_i.
\label{eqn:eR}
\end{align}
\item The error function $\Psi(\bar R,t)$ is bounded by
\begin{equation}\label{eqn:Psi_B}
\frac{1}{4}c_1 \|E_R\|^2 \leq \Psi(\bar R,t) \leq \frac{1}{4}c_2 \|E_R\|^2,
\end{equation}
for some positive constants  $0< c_1\leq c_2$.

\item The attitude error vector $e_R$ satisfies
\begin{equation}\label{eqn:eR_B}
\frac{1}{2}c_1^2(1-\frac{1}{8}\|E_R\|^2) \|E_R\|^2 \leq \|e_R\|^2 \leq \frac{1}{2}c_2^2\|E_R\|^2,
\end{equation}
with the constants $c_1,c_2$ introduced in \eqref{eqn:Psi_B}.

\item Let $Q=R\bar R^T \in\SO$. For any $x,y\in\Re^3$,
\begin{gather}
-x^T(\trs{Q}I_{3\times 3}-Q)x \leq -2 x^T(1-\frac{1}{4}\|E_R\|^2) x, \label{eqn:Qineq0}\\
y^T(\trs{Q}I_{3\times 3} - Q)x  
 \leq 4 \|x\|\|y\|,\label{eqn:Qineq05}\\
(Q-Q^T)^\vee \cdot x\leq \sqrt{2} \|E_R\|\|x\|.\label{eqn:Qineq1}
\end{gather}
\end{enumerate}

\end{lemma}

\subsection{Attitude and Bias Observer}

Consider the following observer presented in~\cite{MahHamITAC08}. 
For positive constants $k_R,k_\gamma$, the attitude and bias observer is formulated as
\begin{align}
\dot{\bar{R}} & =\bar R[\Omega_z-\bar\gamma-k_R e_R]^{\wedge},\label{eqn:R_bar_dot}\\
\dot{\bar\gamma} &= k_\gamma e_R.\label{eqn:gam_bar_dot}
\end{align}
By linearizing the observer dynamics, it has been shown that this observer guarantees local exponential stability when the reference directions are fixed, i.e., $\dot s_i(t)=0$~\cite{MahHamITAC08}. 
Later in~\cite{GriFosITAC12}, the restriction of fixed reference directions is eliminated by introducing a projection operator in the bias observer. 

In this paper, we show that in fact the above observer can handle time-varying reference directions without need for additional modifications.
More specifically, we present an alternative Lyapunov stability analysis without relying on linearization, to show local exponential stability of \eqref{eqn:R_bar_dot} and \eqref{eqn:gam_bar_dot} when the reference direction is time-varying.

\begin{theorem}\label{thm:ASG}
Suppose $k_R > {2d}/{c_1^2}$. 
Then, $(\bar R,\bar\gamma)=(R,\gamma)$ is locally exponentially stable.

More specifically, for $a\in(0,\frac{1}{2})$, choose a constant $\mu$ satisfying
\begin{align}
\mu < \min\big\{ & \sqrt{\frac{c_1}{4k_\gamma}},
    \frac{(1-a)c_1^2k_R-d}{2c_2k_\gamma},\nonumber\\
& \frac{\frac{1}{2}(1-2a)\{(1-a)c_1^2k_R-d\}}{(c_2k_R+\frac{B_\Omega}{2})^2+c_2k_R(1-2a)}\big\}, \label{eqn:mu}
\end{align}
Define the matrices $M_1,M_2,M_3\in\Re^{2\times 2}$ as
\begin{gather*}
    M_1 = \begin{bmatrix} \frac{1}{4}c_1 & -\frac{\sqrt{2}}{2}\mu \\-\frac{\sqrt{2}}{2}\mu & \frac{1}{2k_\gamma}\end{bmatrix},\quad
    M_2 = \begin{bmatrix} \frac{1}{4}c_2 & \frac{\sqrt{2}}{2}\mu \\\frac{\sqrt{2}}{2}\mu & \frac{1}{2k_\gamma}\end{bmatrix},\\
    M_3 ={\footnotesize\selectfont\begin{bmatrix}
            \frac{1}{2}\{(1-a)c_1^2k_R-d\}-\mu c_2 k_\gamma  & -\mu \sqrt{2}(c_2k_R+\frac{B_\Omega}{2}) \\
            -\mu \sqrt{2}(c_2k_R+\frac{B_\Omega}{2}) & 2\mu (1-2a)
    \end{bmatrix}}.
\end{gather*}
Let the constants $\beta$ and $\sigma$ be
\begin{equation}
    \beta=\sqrt{\frac{\lambda_{\max}(M_2)}{\lambda_{\min}(M_1)}},\quad \sigma = \frac{\lambda_{\min}(M_3)}{\lambda_{\max}(M_2)}.\label{eqn:beta_sigma}
\end{equation}
Then, for the trajectories starting from an initial condition satisfying
\begin{equation}
    \Psi(\bar R(0),0) +\frac{1}{2k_\gamma}\| e_\gamma (0)\|^2 \leq 2a c_1,\label{eqn:IC}
\end{equation}
the estimation error exponentially converges as
\begin{equation}\label{eqn:zB}
    \|z(t)\|\leq \beta \|z(0)\| e^{-\frac{\sigma}{2} t},
\end{equation}
where $z=(\|E_R\|,\|e_\gamma\|)\in\Re^2$.
\end{theorem}

Now, we characterize the region of attraction for $(\bar R(0),\bar\gamma(0))\in\SO\times \Re^3$ estimated by \refeqn{IC}. 
It is clear that the region of attraction projected to $\Re^3$ enlarges to $\Re^3$ in the semi-global sense as $k_\gamma\rightarrow\infty$. 
Next, for the region of attraction for the initial attitude estimate, \refeqn{IC} is not useful as $\gamma$, and therefore $e_\gamma(0)$ are not available. 
To remedy this, suppose that the bound of the bias is available as follows.
\begin{assumption}
There is a positive constant $B_\gamma$ satisfying
\begin{equation}\label{eqn:B_gamma}
\|\gamma\| \leq B_\gamma.
\end{equation}
\end{assumption}
In this case, it is reasonable to assume the initial estimate of the bias is selected to satisfy \eqref{eqn:B_gamma}, i.e., $\|\bar\gamma(0)\|\leq B_\gamma$ so that $\|e_\gamma(0)\|\leq 2B_\gamma$.  Assuming that $k_\gamma$ is chosen sufficiently large to satisfy
\begin{equation}
k_\gamma > \frac{B_\gamma^2}{ac_1},\label{eqn:k_gamma_ineq}
\end{equation}
the equation \refeqn{IC} can be rewritten as an inequality for the initial attitude estimate as
\begin{equation}\label{eqn:IC_k_gamma}
\Psi(\bar R(0),0) < 2\parenth{ac_1 - \frac{B_\gamma^2}{k_\gamma}}.
\end{equation}
Applying \refeqn{Psi_B}, a more conservative estimate for the initial attitude estimate guaranteeing exponential convergence is 
\[
\|E_R(0)\|^2  < \frac{8}{c_2}\parenth{ac_1 - \frac{B_\gamma^2}{k_\gamma}}.
\]
For $a\in(0,\frac{1}{2})$, $0<c_1\leq c_2$, and $0<k_\gamma$, the supremum of the right hand side is 4. 
From \eqref{eqn:ER2}, this corresponds to $90^\circ$ of error in the initial attitude estimate. 

In summary, this paper presents the following stability properties of the attitude observer defined by \eqref{eqn:R_bar_dot} and \eqref{eqn:gam_bar_dot}: 
(i) local exponential stability is guaranteed even for time-varying reference directions if the gain $k_R$ is sufficiently large relative to the rate of change of the reference directions; 
(ii) to address the case of time-varying reference directions, there is no need to introduce a projection operator in~\eqref{eqn:gam_bar_dot} as presented in~\cite{GriFosITAC12}; 
(iii) the initial error in the bias estimation can be arbitrarily large provided that $k_\gamma$ is sufficiently large; 
(iv) with regards to the initial attitude estimation error, the inequality \eqref{eqn:IC_k_gamma} guarantees exponential convergence, and it covers at most $90^\circ$ of errors.
All of these are obtained by rigorous Lyapunov stability analysis, and these are unique contributions. 

\section{Numerical Example}

We consider a vehicle equipped with an accelerometer to measure the direction of gravity and a visual sensor to detect feature points. 
The position of the vehicle is available as a function of time as $x(t)=(t,0,0)\in\Re^3$ in the inertial frame. 
The direction of gravity is $e_3=(0,0,1)$, and there are two feature points located at $x_1=(5,0,1)$ and $x_2=(7,-2,0)$. 
Therefore, the three time-varying reference directions are given by
\begin{align*}
    s_1(t) = \frac{x_1 -x(t)}{\|x_1-x(t)\|},\quad
    s_2(t) = \frac{x_2 -x(t)}{\|x_2-x(t)\|},\quad s_3(t)=e_3.
\end{align*}
The weight is chosen as $w_1=w_2=1$ and $w_3=2$, and it can be numerically shown that $c_1=1$, $c_2=4$, and $d=1.14$ for $t\in[0,10]$.

The true attitude trajectory is chosen as $R(t)=\exp(t\hat e_1)\exp(t\hat e_3)\exp(t\hat e_1)$ and $\Omega(t)=(1+\cos t, \sin t - \sin t\cos t, \cos t+\sin^2 t)$ with $R(0)=I_{3\times 3}$. 
The actual gyro bias is $\gamma=(1, 0.5,-1)$. 

\begin{figure}
\centerline{
	\subfigure[Attitude estimation error]{\includegraphics[width=0.7\columnwidth]{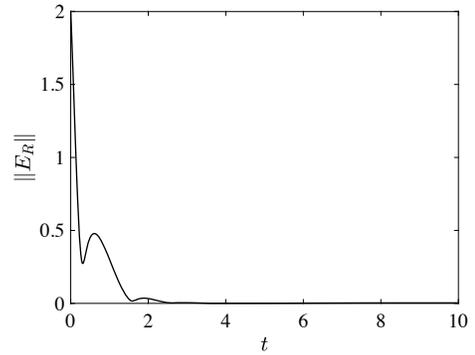}}
}
\centerline{
	\subfigure[Bias estimation error]{\includegraphics[width=0.7\columnwidth]{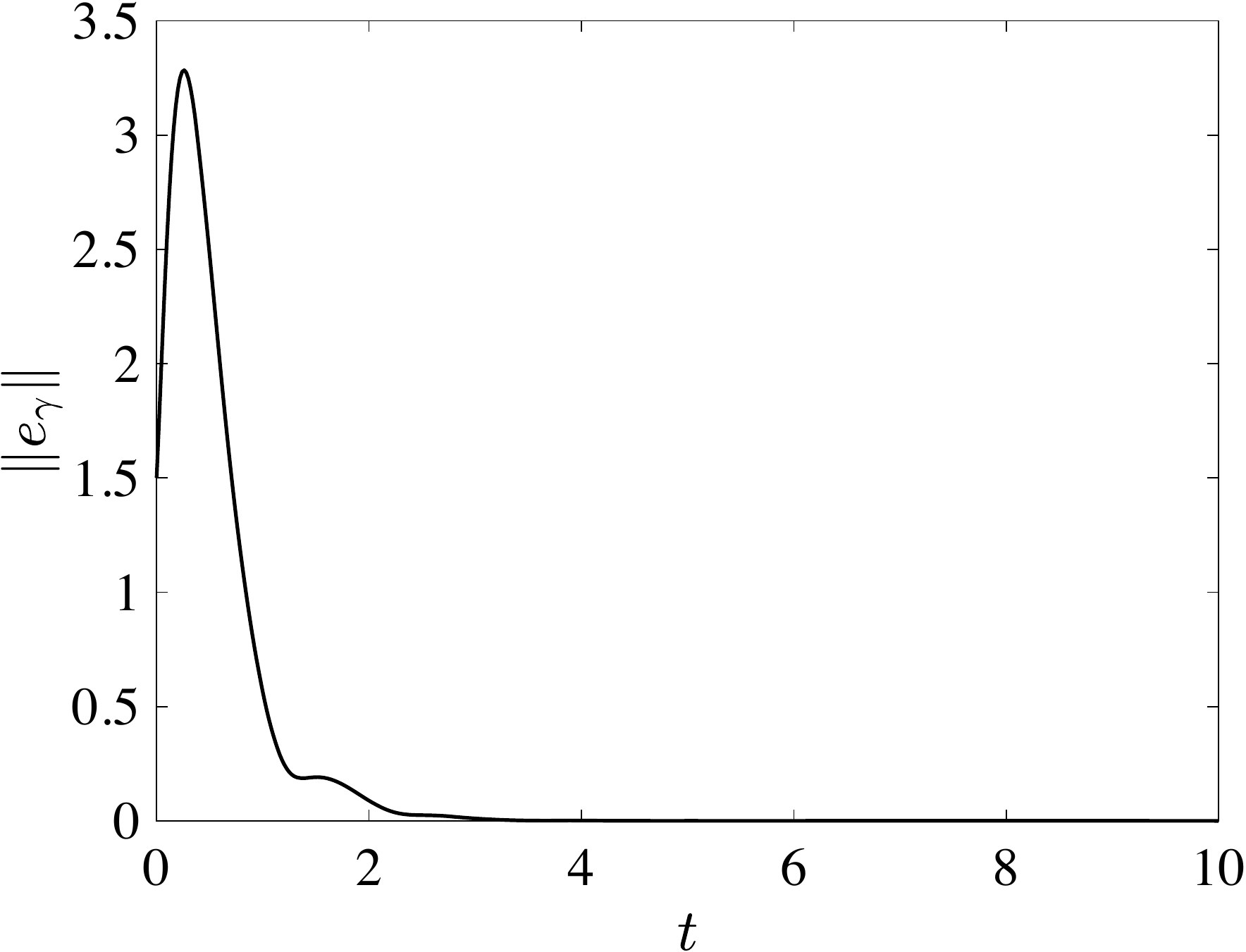}}
}
\caption{Simulation results}\label{fig:NE}
\end{figure}

The initial estimate are $\bar R(0)=\exp(0.5\pi \hat e_1)$ and $\bar \gamma(0)=(0,0,0)$. 
The observer parameters are selected as
\begin{gather*}
\epsilon=0.9,\quad a=0.5\epsilon,\quad B_\gamma=1.65\|\gamma\|,\\
k_R=\frac{2d}{c_1^2\epsilon}=2.53,\quad k_\gamma=\frac{B_\gamma^2}{ac_1\epsilon}=1.65.
\end{gather*}
The corresponding simulation results are presented in Figure \ref{fig:NE}, where the estimation errors converge to zero for time-varying reference directions considered in this example. 

\section{Experimental Results}

The proposed observer is also validated with an indoor attitude estimation experiment.
As it is performed within a building in a crowded urban area, the magnetic field is not consistent. 
As such, the direction measurements are constructed by gravity measured by an acceleration, and visual landmarks captured by a camera. 

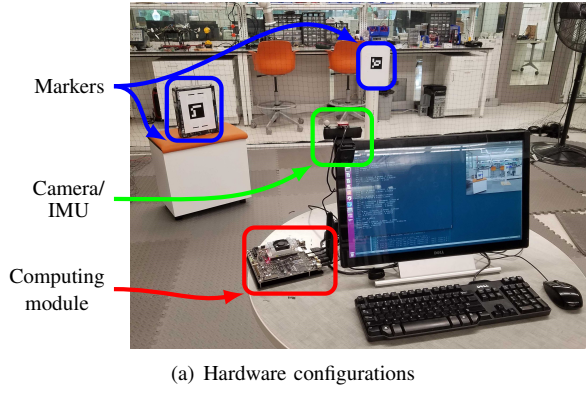
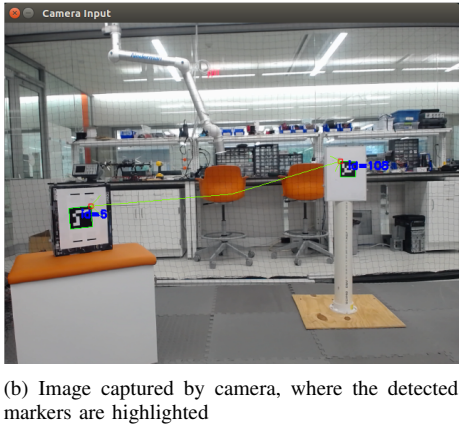
\begin{figure}
    \centerline{
        \subfigure[Hardware configurations]{
            \begin{tikzpicture}
                \node [anchor=east] (marker) at (0.3,3.5) {\footnotesize Markers};
                \node [anchor=east] (imu) at (0.3,2) {\footnotesize \shortstack[c]{Camera/\\IMU}};
                \node [anchor=east] (jetson) at (0.3,0.8) {\footnotesize \shortstack[c]{Computing\\module}};
                \begin{scope}[xshift=0.06\columnwidth]
                    \node[anchor=south west,inner sep=0] (image) at (0,0) {\includegraphics[width=0.7\columnwidth]{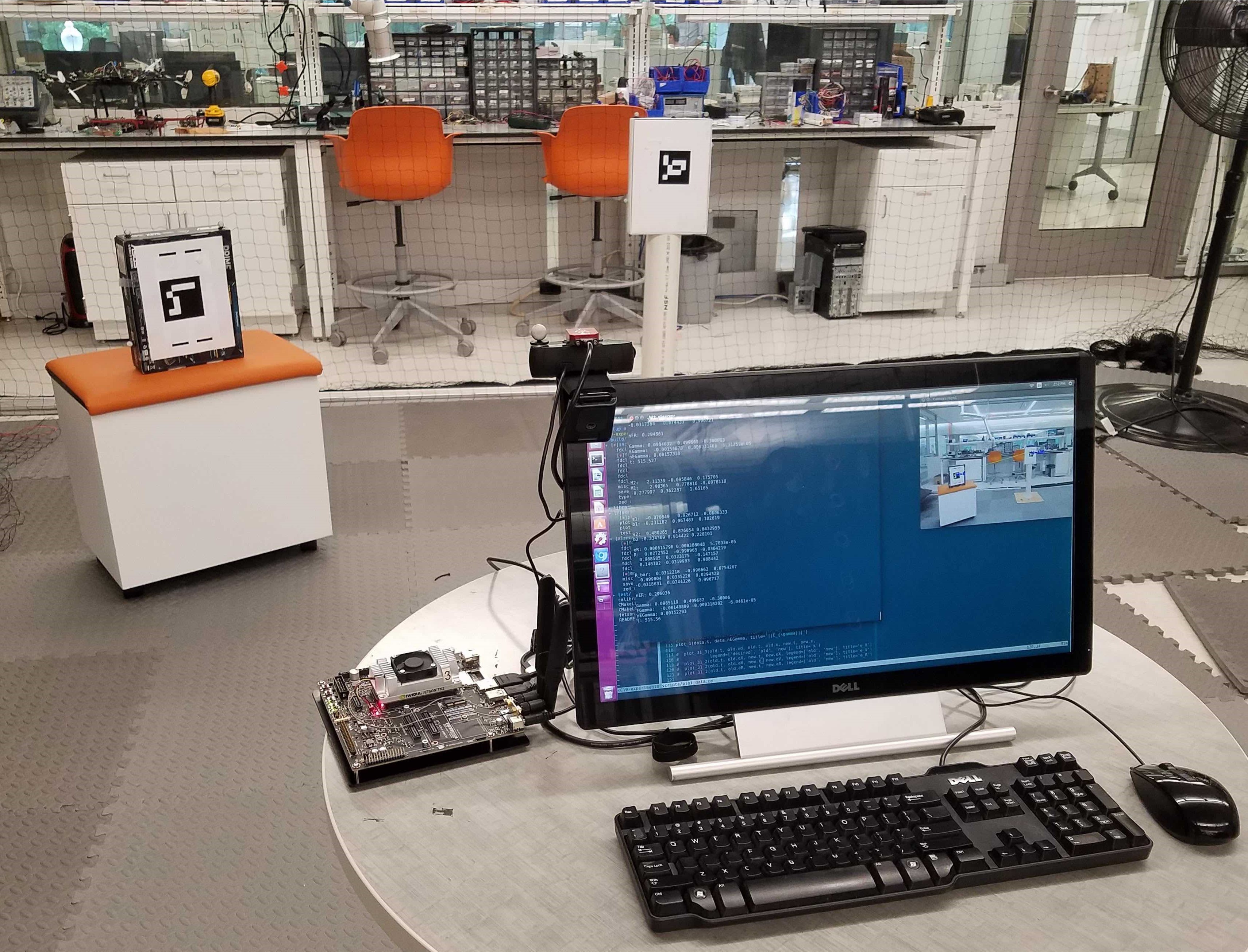}};
                    \begin{scope}[x={(image.south east)},y={(image.north west)}]
                        \draw[blue,ultra thick,rounded corners] (0.08,0.6) rectangle (0.2,0.78);
                        \draw[blue,ultra thick,rounded corners] (0.5,0.75) rectangle (0.58,0.88);
                        \draw [-latex, ultra thick, blue] (marker) to[out=0, in=150] (0.08,0.6);
                        \draw [-latex, ultra thick, blue] (marker) to[out=0, in=150] (0.5,0.88);
                        \draw[green,ultra thick,rounded corners] (0.4,0.53) rectangle (0.53,0.69);
                        \draw [-latex, ultra thick, green] (imu) to[out=0, in=210] (0.4,0.53);
                        \draw[red, ultra thick,rounded corners] (0.25,0.16) rectangle (0.45,0.35);
                        \draw [-latex, ultra thick, red] (jetson) to[out=0, in=190] (0.25,0.16);
                    \end{scope}
                \end{scope}
            \end{tikzpicture}
        }
    }
    \centerline{
        \subfigure[Image captured by camera, where the detected markers are highlighted]{\includegraphics[width=0.7\columnwidth]{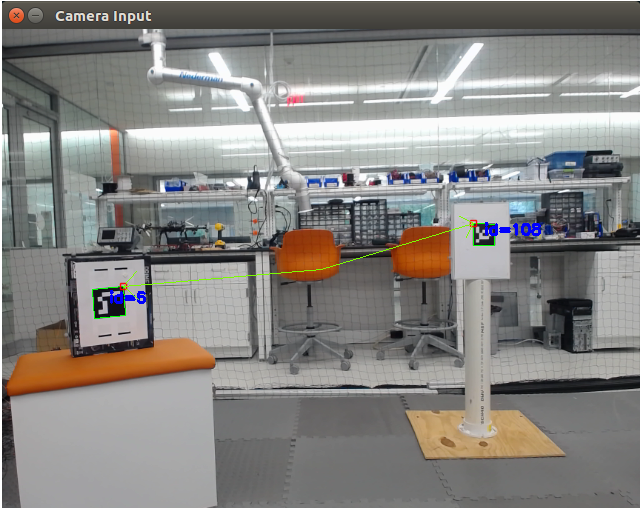}}
    }

    \caption{Experimental setup}
    \label{f:hardware_setup}
\end{figure}

More specifically, the hardware configuration for the presented experiment is as follows.
Two distinctive markers with the pattern of ArUco markers~\cite{Aruco2014} are placed in the lab, and they are fixed in an inertial frame as shown in Figure \ref{f:hardware_setup}a.
A wide angle camera (Logitech C930e 1080P) is used to identify the markers.
It is connected to a computing module (NVidia Jetson TX2) which utilizes the ArUco library and OpenCV~\cite{OpenCV} to detect the markers and compute line of sight represented in the body-fixed frame (refer Figure \ref{f:hardware_setup}b for the image captured by the camera and the detected markers visualized).

Further, a 9-axis inertial measurement unit (VectorNav VN100 IMU) is firmly attached on the camera to measure the direction of the gravitational acceleration and to measure the body angular velocity.
It is connected to Jetson TX2 over a serial port. 
The normalized vector to the each marker from the camera and the normalized direction of the gravitational acceleration are used as $b_1(t)$, $b_2(t)$, and $b_3(t)$, respectively.

The corresponding directions in the inertial frame are measured accurately by a motion capture system composed of five VICON infrared cameras. 
Reflective markers are attached to the camera and each ArUco marker, and their actual orientation and the position is measured by the VICON system at \SI{100}{\hertz}.
The actual position of the camera and each marker is communicated to Jetson TX2 from the VICON server through a wi-fi connection.
This data is used to calculate $s_1(t)$ and $s_2(t)$, while $s_3(t)$ is chosen to be $e_3$.

Finally, the presented observer is implemented at Jetson TX2 via multi-threaded programming in c++, which execute multiple tasks of data acquisition, image processing, observer computation, and data logging simultaneously. 
The observer runs at \SI{30}{\hertz} corresponding to the frame rate of the camera. 
Further, an artificial bias of $\gamma = (0.1, 0.3, -0.2)\,\mathrm{rad/sec}$ is added to the angular velocity measurement to understand the effects of a gyro bias over a short period of time. 
The tests are performed while rotating and translating the camera arbitrarily.

\begin{figure}
    \centerline{
    \subfigure[R: actual (dashed line) and estimated (solid line)]{\includegraphics[width=0.7\columnwidth]{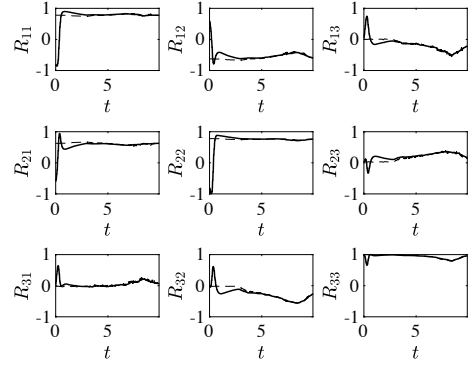}}
    }
    \centerline{
        \subfigure[$\gamma$: actual (dashed line) and estimated (solid line)]{\includegraphics[width=0.7\columnwidth]{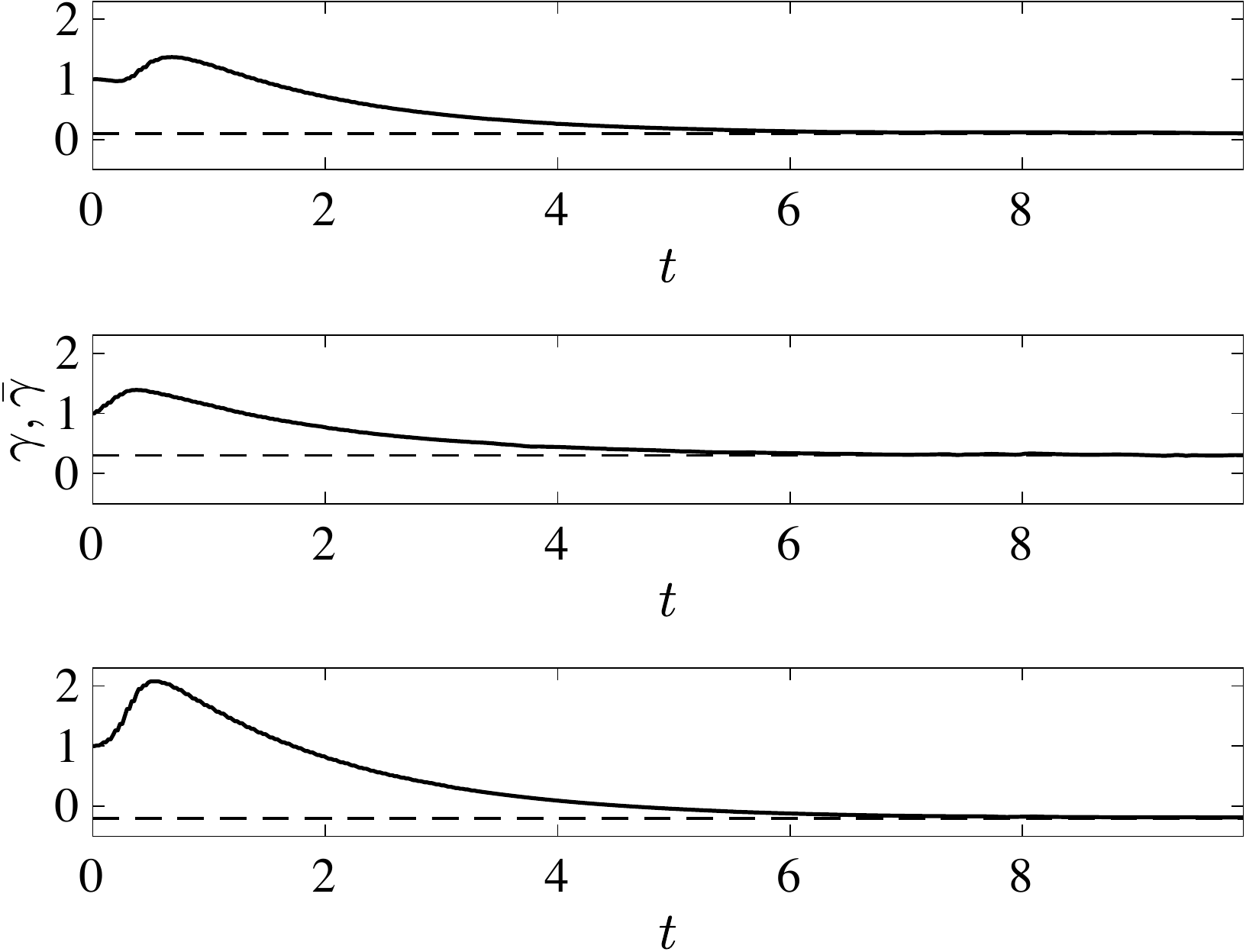}}
    }

\centerline{
	\subfigure[Attitude estimation error]{\includegraphics[width=0.7\columnwidth]{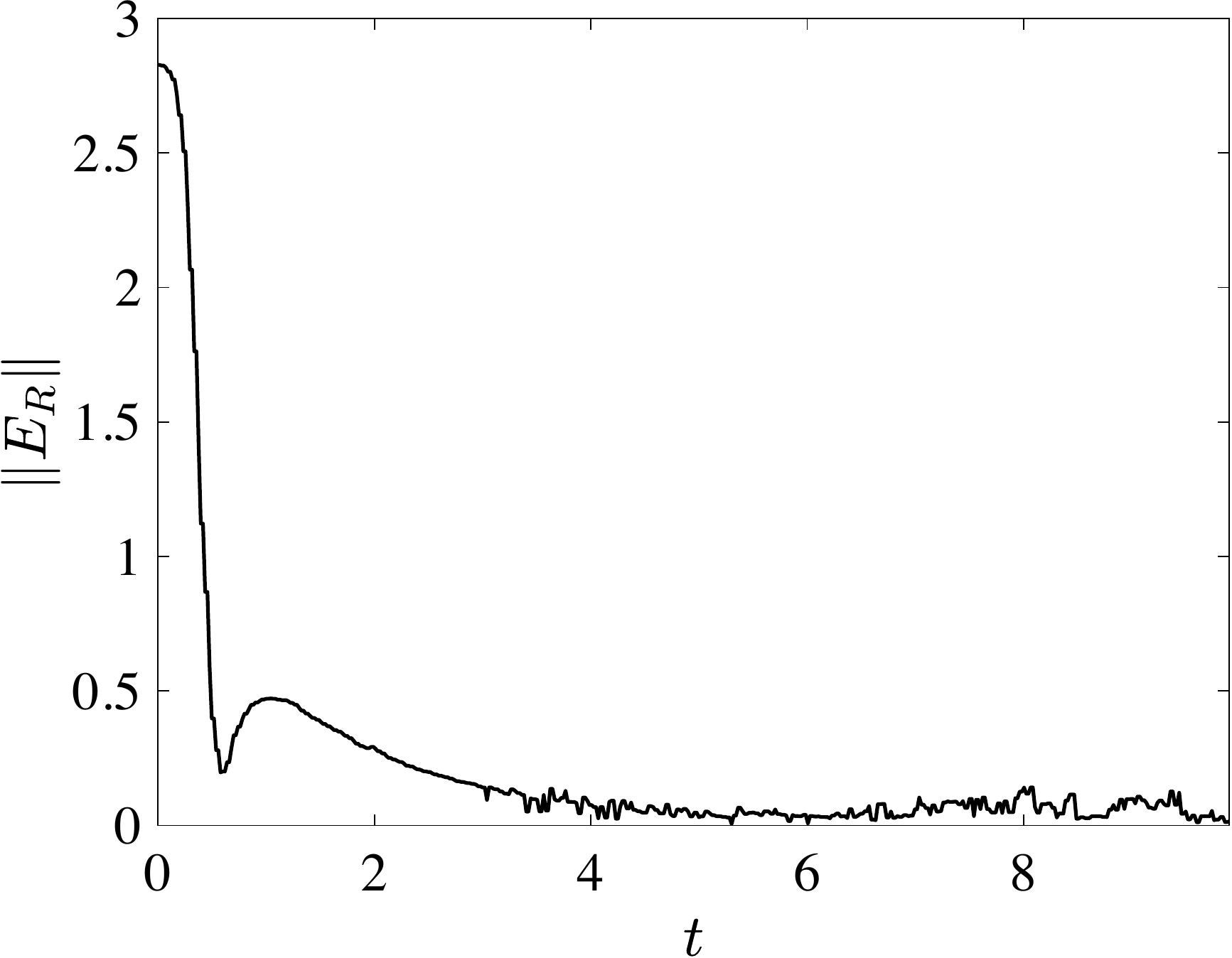}}
}
\centerline{
	\subfigure[Bias estimation error]{\includegraphics[width=0.7\columnwidth]{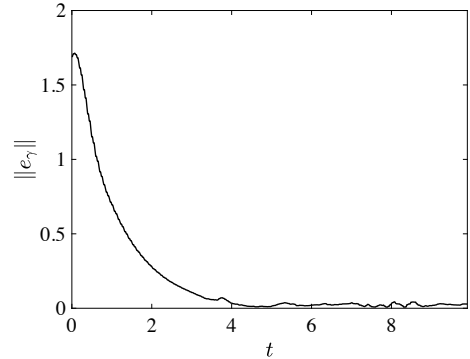}}
}
\caption{Experimental results}\label{fig:NE_exp}
\end{figure}

The experimental results are illustrated in Figure \ref{fig:NE_exp}.
In these plots, the attitude measured by the VICON system is considered as the true attitude, against which the estimation errors are computed. 
For the experimental results presented in this paper, the initial guesses for the observer were set such that $||E_R||=2.828$ and $||E_{\gamma}|| = 1.690$.
The corresponding initial attitude estimation error is close to $180^\circ$, and therefore, it is beyond the presented region of attraction that is possibly conservative. 
However, it shows desirable convergence for both of the attitude estimation error and the gyro bias error. 

\section*{Appendix}

\subsection{Proof of Lemma \ref{lemma:IQ}}
\begin{proof}
The proof is based on the fact that the both sides of the equality is invariant under similarly transformation with rotation matrices. 
Specifically, for any $Q\in\SO$, there exists $P\in\SO$ such that $Q=\exp(\theta \widehat {P e_3})=P\exp(\theta \hat e_3) P^T$, where $e_3=(0,0,1)\in\Sph^2$, i.e., $v_Q$ can be written as $v_Q = P e_3$. 

Substituting this and using \eqref{eqn:R:isometry:0}, it is straightforward to show 
\begin{align*}
    & \langle G(I_{3\times 3}-Q), (I_{3\times 3}-Q)\rangle  \\
&= \langle G(I_{3\times 3}-\exp(\theta \hat e_3)), (I_{3\times 3}-\exp(\theta \hat e_3))\rangle \\
& = 2(G_{11} +G_{22})(1-\cos\theta)\\
&=2(\trs{G} - \langle e_3, Ge_3\rangle)(1-\cos\theta),
\end{align*}
which is equal to the right hand side of \eqref{eqn:IQ} with $v_Q= Pe_3$. 
\end{proof}

\subsection{Proof of Lemma \ref{lemma:important}}
\begin{proof}
The derivative of $\Psi$ with respect to $\bar R$ along $\delta \bar R = \bar R\hat{\bar\eta}$ is given by
\begin{align*}
\D_R \Psi\cdot \delta \bar R  &= \langle GE_R, -\bar R\hat {\bar \eta}\rangle  = -\langle \bar R^TGR, \hat{\bar \eta}\rangle + \langle \bar R^TG\bar R, \hat{\bar \eta}\rangle \\
&= -\frac{1}{2} \langle \bar R^TGR - R^TG\bar R, \hat{\bar \eta}\rangle  = \frac{1}{2} \langle \hat e_R, \hat{\bar \eta}\rangle = e_R \cdot \bar\eta.
\end{align*}
which shows \refeqn{DPsi}. 
The second equality of \refeqn{eR} follows from the identity $(yx^T-xy^T)^\vee=x\times y$ for any $x,y\in\Re^3$. 

Let $Q = R\bar R^T \in \SO$. Then, by Lemma \ref{lemma:IQ}
\begin{align*}
    \Psi(\bar R,t) &= \frac{1}{2}\langle G(t)(R(t)-\bar R), (R(t)-\bar R)\rangle \\
                   &= \frac{1}{2}\langle G(t)(I_{3\times 3} - Q), (I_{3\times 3} - Q)\rangle \\
                   &= (\trs{G(t)} - \langle v_Q, G(t)v_Q\rangle) (1-\cos\theta),
\end{align*}
where $\theta \in \mathbb R$ and ${v_Q} \in \Sph^2$ such that $Q = \exp(\theta \widehat {v_Q})$.  Since 
\begin{equation}\label{eqn:ER2}
\|E_R\|^2 = \|R-\bar R\|^2 = \|I_{3\times 3} - Q\|^2 = 4(1-\cos\theta),
\end{equation}
we have
\begin{align*}
    \Psi(\bar R,t) & = \frac{1}{4} (\trs{G(t)} - \langle v_Q, G(t)v_Q\rangle) \| E_R\|^2.
\end{align*}
Let $\lambda_1(t), \lambda_2(t), \lambda_3(t)\in\Re$ be the eigenvalues of $G(t)$ in ascending order, i.e., $\lambda_1(t)=\lambda_{\min}(G(t))$ and $\lambda_2(t)=\lambda_{\max}(G(t))$.
As $\trs{G(t)}=\sum_{i=1}^3 \lambda_i(t)$ and $\|v_Q\|=1$, we have
\begin{align*}
    \lambda_1(t)+\lambda_2(t) \leq \trs{G(t)} - \langle v_Q, G(t)v_Q\rangle \leq \lambda_2(t)+\lambda_3(t).
\end{align*}
From \eqref{eqn:lambda_G}, $\lambda_1(t)+\lambda_2(t)\geq c \triangleq c_1$. 
Also, it is well known that the spectral radius of a matrix is bounded by any matrix norm. 
The Froubenius norm of $G(t)$ is bounded as $G(t)$ is composed of the finite elements of $s_i\in\Sph^2$. 
Therefore, there exists a positive constant $c_2$ such that $\lambda_2(t)+\lambda_3(t) \leq c_2$ always. 
These show \eqref{eqn:Psi_B}.

Let $Q = R\bar R^T \in \SO$. Then,
\begin{align}
\|e_R\|^2 
&= \frac{1}{2}\| R^TG\bar R - \bar R^T G R\|^2\nonumber \\
&= \frac{1}{2}\|GQ^T- QG\|^2 \label{here:eR:0} \\
&=\frac{1}{2}\| (I_{3\times3}-Q)G - G(I_{3\times3}-Q^T)\|^2 \label{here:eR}\\
&= \|G(I_{3\times3}-Q^T)\|^2 - \langle (I_{3\times 3}-Q)G, G(I_{3\times3}-Q^T)\rangle. \nonumber
\end{align}
Since the inequalities \eqref{eqn:eR_B} are invariant under similarity transformation by $\SO$ matrices, similar with the proof of Lemma \ref{lemma:IQ}, we may assume that $Q=\exp(\theta\hat e_3)$ for some $\theta \in [0,2\pi)$. 
Then, by Lemma \ref{lemma:IQ} for $G^2$,
\begin{align}\label{IGQ}
\|G&(I_{3\times 3}-Q^T)\|^2  \\
&= \langle G(I_{3\times3}-Q^T), G(I_{3\times3}-Q^T)\rangle \nonumber\\
&= \langle G^2(I_{3\times3}-Q^T),I_{3\times3} - Q^T\rangle \nonumber \\
&= 2((G_{11})^2  + (G_{13})^2 + (G_{22})^2 + (G_{23})^2)(1-\cos\theta ). \nonumber
\end{align}
It is straightforward to compute
\begin{align*}
&\langle (I_{3\times3}-Q)G, G(I_{3\times3}-Q^T)\rangle \\
&= ((G_{11})^2 + (G_{22})^2)(1-\cos\theta)^2 - 2G_{11}G_{22}(1-\cos^2\theta).
\end{align*}
Putting these all together with \eqref{eqn:ER2}, we get
\begin{align}
2\frac{\| e_R\|^2}{\|  E_R\|^2} &= (G_{13})^2 + (G_{23})^2  \nonumber\\
&\quad+ (G_{11} + G_{22})^2\left (   1-\frac{1}{8} \| E_R\|^2 \right). \label{eR:equality}
\end{align}
It follows
\begin{align*}
2\frac{\| e_R\|^2}{\|  E_R\|^2} &\geq (G_{11} + G_{22})^2\left (   1-\frac{1}{8} \| E_R\|^2 \right)\\
&\geq (\trs{G} - \lambda_{\rm max}(G))^2 \left (   1-\frac{1}{8} \| E_R\|^2 \right),
\end{align*}
which implies the left inequality in \eqref{eqn:eR_B}. From \eqref{eR:equality},
\begin{align*}
2\frac{\|e_R\|^2}{\|E_R\|^2} &\leq (G_{13})^2 + (G_{23})^2  + (G_{11}+G_{22})^2.
\end{align*}
Let $F=\trs{G}I_{3\times 3} - G\in\Re^3$. 
Then, it can be shown that the right hand side of the above inequality is equal to $e_3^T F^2 e_3$, which is less than or equal to $\lambda_{\max}(F^2)=\lambda_{\max}^2(F)=(\trs{G}-\lambda_{\min} (G))^2$. 
This shows the right inequality in \eqref{eqn:eR_B}.

Let us now show (iv). 
The thee inequalities \eqref{eqn:Qineq0} -- \eqref{eqn:Qineq1} are invariant under similarity transformation by $\SO$ matrices. 
Hence, we may assume that $Q$ is of the form $Q=\exp(\theta\hat e_3)$ without loss of generality. 
Then, \eqref{eqn:ER2} implies $\trs{Q} = 1-2\cos\theta = 3-\frac{1}{2}\|E_R\|^2$ from which \eqref{eqn:Qineq0} follows. 
The inequality \eqref{eqn:Qineq05} follows from the triangle inequality and the Cauchy-Schwarz inequality. 
Finally, it is easy to compute
\[
(Q-Q^T)^\vee \cdot x = -2 x_2\sin \theta \leq 2 |\sin\theta | \|x\|.
\]
Since $2 \cos\theta + \sin^2\theta \leq 2$ for all $\theta \in [0,2\pi]$, it follows that $2|\sin\theta | \leq \sqrt 2\sqrt{4(1-\cos\theta )} = \sqrt 2\|E_R\|$ from which \eqref{eqn:Qineq1} follows.
\end{proof}

\subsection{Proof of Theorem \ref{thm:ASG}}
\begin{proof}
Let a Lyapunov function be
\begin{equation}
V_0(\bar R,\bar\gamma,t) = \Psi(\bar R,t) + \frac{1}{2k_\gamma} \| e_\gamma \|^2.\label{eqn:V0}
\end{equation}
From \refeqn{Psi_B}, this is positive-definite and decrescent about $z=0$. 
For $a\in(0,\frac{1}{2})$, define an open domain about $z=0$ as
\begin{equation}
    D=\{ (\bar R,\bar\gamma)\in\SO\times\Re^3\,|\, V_0(\bar R, \bar\gamma,t) < 2 a c_1 \}.\label{eqn:D}
\end{equation}

Using \refeqn{DPsi}, and \refeqn{R_bar_dot}, the time-derivative of $V_0$ is given by
\begin{align*}
    \dot V_0(\bar R,\bar\gamma,t) & = e_R\cdot(-\Omega + \Omega_z -\bar\gamma -k_R e_R) +\frac{1}{k_\gamma} e_\gamma \cdot\dot e_\gamma\\
                                  & \quad + \frac{1}{2}\pair{ \dot G E_R, E_R}.
\end{align*}
Substituting \refeqn{Wz} and \refeqn{gam_bar_dot} with $\dot{ e}_\gamma =-\dot{\bar\gamma}$, and from \eqref{eqn:d},
\begin{equation}
    \dot V_0(\bar R,\bar\gamma,t) \leq  -k_R \|e_R\|^2 + \frac{1}{2}d \|E_R\|^2.\label{eqn:V0_dot}
\end{equation}
Using \refeqn{Psi_B}, 
\begin{equation}\label{eqn:ER_IC}
\|E_R(t)\|^2 \leq \frac{4}{c_1}\Psi(t) \leq 8a,
\end{equation}
where the second inequality is obtained by \eqref{eqn:D}. Thus, from \refeqn{eR_B},
\begin{gather}
\frac{1}{2}c_1^2(1-a)\|E_R\|^2\leq \|e_R\|^2 \leq 4a c_2^2 < 4c_2^2.\label{eqn:eR_IC}
\end{gather}
Substituting this into \eqref{eqn:V0_dot}
\begin{gather}
    \dot V_0(\bar R,\bar\gamma,t) \leq -\frac{1}{2}\{c_1^2(1-a)k_R-d\}\|E_R\|^2.\label{eqn:V0_dot_IC}
\end{gather}
Therefore, if $k_R > \frac{2d}{c_1^2 }> \frac{d}{c_1^2 (1-a)}$, $\dot V_0$ is negative-semidefinite, which follows that $z=0$ is stable in the sense of Lyapunov. 
Also, the given domain $D$ is positively invariant. i.e., if $(\bar R(0),\bar\gamma(0))\in D$ then $(\bar R(t),\bar \gamma(t))\in D$ for all $t\geq 0$.

For $\mu>0$, let the augmented Lyapunov function be 
\begin{equation}
V(\bar R,\bar\gamma,t) = V_0(\bar R,\bar\gamma,t)+ \mu (R  e_\gamma ) \cdot (Q-Q^T)^\vee.
\label{eqn:V}
\end{equation}
Using \refeqn{Psi_B} and \refeqn{Qineq1}, it satisfies
\[
z^T M_1 z \leq V(\bar R,\bar\gamma,t) \leq z^T M_2 z,
\]
where the matrices $M_1$ and $M_2$ are positive-definite due to \refeqn{mu}, and therefore, $V(\bar R,\bar\gamma,t)$ is positive-definite and decrescent. 

Next, we derive the time-derivative of $V(\bar R,\bar\gamma,t)$. 
Utilizing \refeqn{V0_dot_IC}, we just need to find the time-derivative of the last, cross term in \refeqn{V}. First, using \refeqn{R_dot} and \refeqn{R_bar_dot}, 
\begin{align*}
\dot Q & = R\hat\Omega \bar R ^T - R(\Omega_z-\bar\gamma-k_R e_R)^\wedge \bar R^T\\
& = R (- e_\gamma  + k_R e_R)^\wedge \bar R^T\\
& = \{R(- e_\gamma  + k_R e_R)\}^\wedge R\bar R^T
\triangleq \hat\omega_Q Q,
\end{align*}
with $\omega_Q = R(- e_\gamma  + k_R e_R)\in\Re^3$. Thus,
\begin{align}
\frac{d}{dt}\{( R  e_\gamma ) &\cdot (Q-Q^T)^\vee\}
 =R(\hat\Omega e_\gamma  - k_\gamma e_R)\cdot (Q-Q^T)^\vee\nonumber\\
&\quad + (R e_\gamma )\cdot (\hat\omega_Q Q +Q^T\hat\omega_Q)^\vee.\label{eqn:cross_dot}
\end{align}
From \refeqn{Qineq1}, \refeqn{BW} and \refeqn{eR_B}, the first term of the right hand side of \refeqn{cross_dot} satisfies
\begin{align}
R(\hat\Omega e_\gamma  & - k_\gamma e_R)\cdot (Q-Q^T)^\vee\nonumber\\
&\leq \sqrt{2}(B_\Omega\| e_\gamma \| +\frac{k_\gamma}{\sqrt 2} c_2\|E_R\| ) \|E_R\|.
\label{eqn:cross_dot_1}
\end{align}
Using the hat map identity, $(\hat x A + A^T\hat x)^\vee = (\trs{A}I_{3\times 3}-A)x$ for any $x\in\Re^3$ and $A\in\Re^{3\times 3}$,  the last term of \refeqn{cross_dot} is expanded as
\begin{align*}
(R e_\gamma ) &\cdot (\hat\omega_Q Q +Q^T\hat\omega_Q)^\vee
= (R e_\gamma )^T(\trs{Q}I_{3\times 3}-Q)\omega_Q\\
&=(R e_\gamma )^T(\trs{Q}I_{3\times 3}-Q)(R(- e_\gamma  + k_R e_R)).
\end{align*}
In view of \refeqn{Qineq0} and \refeqn{Qineq05},
\begin{align*}
(R e_\gamma ) &\cdot (\hat\omega_Q Q +Q^T\hat\omega_Q)^\vee\\
&\leq -2(1-\frac{1}{4}\|E_R\|^2)\| e_\gamma \|^2 + 4k_R \| e_\gamma \|\|e_R\|.
\end{align*}
Further from \refeqn{eR_B} and \refeqn{ER_IC},
\begin{align}
(R e_\gamma ) &\cdot (\hat\omega_Q Q +Q^T\hat\omega_Q)^\vee\nonumber\\
&\leq -2(1-2a)\| e_\gamma \|^2 + 2\sqrt{2} k_Rc_2 \| e_\gamma \|\|E_R\|.
\label{eqn:cross_dot_2}
\end{align}
With \refeqn{cross_dot_1} and \refeqn{cross_dot_2}, we obtain an upper bound of \refeqn{cross_dot}, which is combined with \refeqn{V0_dot_IC} to obtain
\begin{equation}
\dot V(\bar R,\bar\gamma,t)\leq - z^T M_3 z\leq -\sigma V(\bar R,\bar\gamma,t),
\end{equation}
where the matrix $M_3\in\Re^{2\times 2}$ is positive-definite from \refeqn{mu}. Therefore, $z=0$ is exponentially stable, and it is straightforward to show
\[
V(\bar R(t),\bar\gamma(t),t) \leq V(\bar R(0),\bar\gamma(0),0) \exp(-\sigma t),
\]
which implies \refeqn{zB}.
\end{proof}



\end{document}